# NONSTABLE K-THEORY FOR $\mathcal{Z}$-STABLE $C^*$-ALGEBRAS

XINHUI JIANG

ABSTRACT. Let $\mathcal{Z}$ denote the simple limit of prime dimension drop algebras that has a unique tracial state (cf. Jiang and Su [11]). Let $A \neq 0$ be a unital $C^*$-algebra with $A \cong A \otimes \mathcal{Z}$. Then the homotopy groups of the group $\mathcal{U}(A)$ of unitaries in $A$ are stable invariants, namely, $\pi_i(\mathcal{U}(A)) \cong K_{i-1}(A)$ for all integer $i \geq 0$. Furthermore, $A$ has cancellation for full projections, and satisfies the comparability question for full projections. Analogous results hold for non-unital $\mathcal{Z}$-stable $C^*$-algebras.

## 0. Introduction and summary of results

Let $\mathcal{Z}$ denote the only simple limit of prime dimension drop algebras that has a unique tracial state (cf. Jiang and Su [11]). A $C^*$-algebra $A$ is called $\mathcal{Z}$-stable, if $A \cong A \otimes \mathcal{Z}$. In this note, we study nonstable K-theory for $\mathcal{Z}$-stable $C^*$-algebras. Our motivations first came from works on approximately divisible $C^*$-algebras (cf. Blackadar, Kumjian and Rordam [4] and the references therein).

At first glance, these two classes of algebras might seem quite unrelated. For example, the algebra $\mathcal{Z}$ itself is known to be $\mathcal{Z}$-stable ([11]). Since $\mathcal{Z}$ has no non-trivial projections, it is certainly *not* approximately divisible. It is quite easy to check that some properties (related to projections) of approximately divisible algebras (cf. Theorem 1.4(b), (e), (f) in [4]) fail on $\mathcal{Z}$.

On the other hand, some very interesting approximately divisible $C^*$-algebras are known to be $\mathcal{Z}$-stable (cf. Theorem 5 of [11]). In fact, to this date we know of no obstructions to $\mathcal{Z}$-stability for approximately divisible (nuclear) $C^*$-algebras. Furthermore, non-zero $\mathcal{Z}$-stable $C^*$-algebras are, in a certain sense (cf. Remark 1.4), "fibrewise approximately divisible", and certain properties of approximately divisible $C^*$-algebras do persist in the class of $\mathcal{Z}$-stable algebras (cf. Theorem 3 of Gong, Jiang and Su [9]). In this note, we establish the following results (compare Theorem 1.4(d), Propositions 3.10, and Proposition 3.11 of [4]):

**Theorem 1.** *Let $A$ be a unital $\mathcal{Z}$-stable $C^*$-algebra. Then:*

(a) *$A$ has cancellation for full projections: If $p$ and $q$ are two full projections in $A$, and $[p] = [q] \in K_0(A)$, then $e$ and $f$ are (Murray-von Neumann) equivalent.*

(b) *$A$ satisfies the comparability question for full projections: If $p$ and $q$ are two full projections in $A$, and $\tau(p) < \tau(q)$ for all quasi-traces $\tau$ on $A$, then $p$ is equivalent to a subprojection of $q$.*

**Theorem 2.** *Let $A$ be a unital $\mathcal{Z}$-stable $C^*$-algebra. Then the natural map $\mu : \mathcal{U}(A)/\mathcal{U}_0(A) \to K_1(A)$ is an isomorphism, where $\mathcal{U}(A)$ denotes the unitary group of $A$ and $\mathcal{U}_0(A)$ the connected component of the identity.*

Our proof of Theorem 1 will rely upon Theorem 2 (to prove Theorem 1.(b), we also need results from Blackadar and Rordam [5]).

Naturally, Theorem 2 leads us to the higher homotopy groups of $\mathcal{U}(A)$. This was also inspired by the results in Zhang [17], where an earlier (and different) notion of approximate divisibility was proposed. We prove that these groups are also stable invariants:

**Theorem 3.** *Let $A \neq 0$ be a unital $\mathcal{Z}$-stable $C^*$-algebra. Then for any integer $i \geq 0$,*

$$\pi_i(\mathcal{U}(A)) = \begin{cases} K_0(A) & \text{if } i \text{ is odd,} \\ K_1(A) & \text{if } i \text{ is even.} \end{cases}$$

Following [17], we also calculate the homotopy groups of the Grassmann space of $A$ when $A$ is unital simple and $\mathcal{Z}$-stable. See Proposition 3.8 for details.

Results similar to Theorem 3 have been established earlier for various classes of $C^*$-algebras (see, for example, Cuntz [7], Rieffel [12], Thomsen [15], and Zhang [17], [18]). We refer the reader to [15] and [17] for careful discussions of earlier works. Our approach is based, as in [15], on a fundamental observation in K-theory (see (1.6) below), and is remarkably elementary (see the end of §1 for a outline). On the other hand, our result covers many, if not most, of the (unital) separable $C^*$-algebras treated earlier. This, we think, is a clear indication of the importance of $\mathcal{Z}$-stability.

To save notation, we treat mainly unital $C^*$-algebras. Some analogous results hold also in the non-unital case. These will be discussed in remarks.

*Throughout this note, any $C^*$-algebra $A$ is assumed to be non-trivial, that is, $A \neq 0$.*

The rest of this note is organized as follows: In §1, we establish notation, recall some basic facts about the algebra $\mathcal{Z}$, and outline the strategy to be followed in §2. The main technical part of this note is contained in §2, where we construct and study an interesting map and prove Theorem 3 (hence also Theorem 2). An extension to the non-unital case is also discussed there. Finally, in §3, we apply Theorem 2 to study full projections, and prove Theorem 1.

We are grateful to George Elliott, Guihua Gong and Hongbing Su for stimulating conversations.

## 1. Preliminaries

In this section, we first establish some notation to be used in this note and recall the most basic facts about the algebra $\mathcal{Z}$. Then we comment on possible relations between $\mathcal{Z}$-stability and approximate divisibility, and outline the basic strategy of this note.

**Notation 1.1.** Let $A$ be a $C^*$-algebra and $m$ and $n$ two positive integers. The following notation will be used in this note:

1. $\mathbf{M}_n$ denotes the algebra of all $n \times n$ complex matrices, with identity element $\mathbf{1}_n$, zero element $\mathbf{0}_n$, and the standard matrix unit $\{\mathbf{e}_{i,j} : 1 \leq i, j \leq n\}$.

2. $\mathbf{M}_n(A)$ denotes the algebra of all $n \times n$ matrices with entries in $A$. We identify $A \otimes \mathbf{M}_n$ with $\mathbf{M}_n(A)$ in the usual way. Namely, we identify $\sum_{i,j} a_{i,j} \otimes \mathbf{e}_{i,j}$ in $A \otimes \mathbf{M}_n$ with (the $n \times n$ matrix) $[a_{i,j}]$ in $\mathbf{M}_n(A)$. In particular, for any $a \in A$,

we have:
$$a \otimes \mathbf{1}_n = \begin{bmatrix} a & 0 & \cdots & 0 \\ 0 & a & \cdots & 0 \\ \vdots & \vdots & & \vdots \\ 0 & 0 & \cdots & a \end{bmatrix} \in \mathbf{M}_n(A).$$

Also, if $A$ is unital, we identify $\mathbf{M}_k$ canonically with the unital subalgebra $1_A \otimes \mathbf{M}_k$ of $\mathbf{M}_k(A)$.

3. For any $a \in \mathbf{M}_m(A)$ and $b \in \mathbf{M}_n(A)$, we write:
$$a \oplus b = \text{diag}(a, b) = \begin{bmatrix} a & 0 \\ 0 & b \end{bmatrix} \in \mathbf{M}_{m+n}(A).$$

4. Let $Z_{m,n}(A)$ denote the $C^*$-algebra of all continuous functions $f : [0,1] \to \mathbf{M}_{mn}(A)$ with
$$f(0) = a \otimes \mathbf{1}_n, \quad \text{for some} \quad a \in \mathbf{M}_m(A),$$
and,
$$f(1) = b \otimes \mathbf{1}_m, \quad \text{for some} \quad b \in \mathbf{M}_n(A).$$

If $A = \mathbb{C}$, we shall denote $Z_{m,n}(A)$ simply by $Z_{m,n}$. (The same algebra is denoted by $\mathbf{I}[m, mn, n]$ in [11].) It is easy to see that
$$Z_{m,n}(A) \cong A \otimes Z_{m,n}.$$
From now on, we shall use these two notations interchangeably.

The algebra $Z_{m,n}$ is called a prime dimension drop algebra, if $m$ and $n$ are relatively prime. The algebra $\mathcal{Z}$ constructed in [11] is the only simple $C^*$-algebra with a unique tracial state which is a limit of prime dimension drop algebras with unital connecting maps. Since each prime dimension drop algebra is nuclear, so is $\mathcal{Z}$.

The following two properties of $\mathcal{Z}$ will be basic to the discussions in this note. They can be found in [11].

**Proposition 1.2.** *([11]) Let $Z_{m,n}$ be a prime dimension drop algebra. Then there exists a unital injective $*$-homomorphism from $Z_{m,n}$ into $\mathcal{Z}$.*

This follows from the construction of $\mathcal{Z}$ (cf. proof of Proposition 2.7 in [11]) and the uniqueness theorem (Theorem 6.2) of [11].

**Theorem 1.3.** *(Theorem 4 of [11]) $\mathcal{Z} \otimes \mathcal{Z} \cong \mathcal{Z}$ and $\mathcal{Z}^{\otimes \infty} \cong \mathcal{Z}$.*

Recall that $\mathcal{Z}^{\otimes \infty}$ is the limit of the sequence $(\mathcal{Z}^{\otimes n}, \iota_n)$, where $\mathcal{Z}^{\otimes n} = \mathcal{Z} \otimes \cdots \otimes \mathcal{Z}$ ($n$ copies of $\mathcal{Z}$) and $\iota_n$ is the canonical embedding map: $\iota_n(a) = a \otimes 1_\mathcal{Z}$ for $a \in \mathcal{Z}^{\otimes n}$.

*Remark 1.4.* It follows from Theorem 1.3 that, if a $C^*$-algebra $A$ is $\mathcal{Z}$-stable, there is an isomorphism $\theta : A \otimes \mathcal{Z}^{\otimes \infty} \cong A$. Let $A_n = \theta(A \otimes \mathcal{Z}^{\otimes n})$ for each $n \geq 0$. Then it is easy to check that:

1. $\{A_n\}$ is an increasing sequence of unital $C^*$-subalgebras of $A$, and $\cup_n A_n$ is dense in $A$, and

2. for each $n$, there is an isomorphism $\theta_n : A_n \otimes \mathcal{Z} \cong A_{n+1}$ such that the following diagram commutes:

$$\begin{array}{ccc} A_n & \xrightarrow{\cdot \otimes 1_{\mathcal{Z}}} & A_n \otimes \mathcal{Z} \\ \downarrow \text{incl.} & & \downarrow \theta_n \\ A_{n+1} & =\!\!=\!\!= & A_{n+1}. \end{array} \quad (1.1)$$

These are reminiscent of the structure of approximately divisible $C^*$-algebras. By Theorem 1.3 of [4], a unital separable $C^*$-algebra $A$ is approximately divisible if, and only if, there is an increasing sequence $\{A_n\}$ of unital $C^*$-subalgebras of $A$ with the following properties:

- $\cup_n A_n$ is dense in $A$, and
- for each $n$, there is an finite dimensional $C^*$-algebra $F_n = \oplus_i \mathbf{M}_{k_i}$, with $k_i \geq 2$ for all $i$, and a unital morphism $\theta_n : A_n \otimes F_n \to A_{n+1}$ such that the following diagram commutes:

$$\begin{array}{ccc} A_n & \xrightarrow{\cdot \otimes 1_{F_n}} & A_n \otimes F_n \\ \downarrow \text{incl.} & & \downarrow \theta_n \\ A_{n+1} & =\!\!=\!\!= & A_{n+1}. \end{array} \quad (1.2)$$

Compressing the sequence if necessary, we can require that the sizes $k_i$ in $F_n$ be greater than any prescribed integer. This property has been crucial in analyzing the structure of a large class of interesting $C^*$-algebras, see [2], [13], [14], and [4] for details. See also [17] for another notion of approximate divisibility.

Note that a prime dimension drop algebra $Z_{m,n}$ is a continuous field of $C^*$-algebras over $[0,1]$, where each fibre is a full matrix algebra $\mathbf{M}_k$ with $k = m$, $n$, or $mn$. Therefore, by Proposition 1.2, the algebra $\mathcal{Z}$ contains continuous fields of full matrix algebras with sizes greater than any prescribed integer. In this sense, $\mathcal{Z}$-stable algebras are "fibrewise approximately divisible".

The following question arises naturally from our discussion.

**Question 1.5.** If a $C^*$-algebra $A$ is approximately divisible, is it $\mathcal{Z}$-stable?

Note that this question is closely related to Elliott's classification program for nuclear $C^*$-algebras ([8]). If $A$ is simple, nuclear, and approximately divisible, then it follows from Theorem 1.4 of [9] that $K_0(A)$ is weakly unperforated. Hence, by Lemmas 2.11 and 2.12 of [11] and Theorem 1 of [9], $A$ and $A \otimes \mathcal{Z}$ have the same Elliott invariants.

In the rest of this section, we outline our approach to Theorem 3. First, we explain the notation.

*Remark 1.6.* For any unital $C^*$-algebra $A$, in this note we consider its unitary group $\mathcal{U}(A)$ as a pointed topological space, with $1_A$ as the basepoint. In particular, homotopy groups are calculated with this basepoint.

Accordingly, if $A$ and $B$ are two unital $C^*$-algebras, then by a map $\psi : \mathcal{U}(A) \to \mathcal{U}(B)$, we mean a continuous and basepoint-preserving map. The same restriction also applies to any homotopy between two such maps.

Next, we recall a basic fact in K-theory. Let $A$ be any unital $C^*$-algebra. Let

$$W_t = \begin{bmatrix} \cos(\pi t/2) & \sin(\pi t/2) \\ -\sin(\pi t/2) & \cos(\pi t/2) \end{bmatrix}, \quad \forall t \in [0,1]; \tag{1.3}$$

and,

$$\rho_t(u) = \begin{bmatrix} u & 0 \\ 0 & 1 \end{bmatrix} W_t^* \begin{bmatrix} 1 & 0 \\ 0 & u \end{bmatrix} W_t, \quad \forall u \in \mathcal{U}(A). \tag{1.4}$$

Then, for any $u$, $\rho(u)$ is a homotopy between the following two unitaries:

$$\rho_0(u) = \begin{bmatrix} u & 0 \\ 0 & u \end{bmatrix}, \quad \text{and} \quad \rho_1(u) = \begin{bmatrix} u^2 & 0 \\ 0 & 1 \end{bmatrix}. \tag{1.5}$$

More generally, for any integer $k \geq 1$, one can construct a homotopy $\rho : [0,1] \times \mathcal{U}(A) \to \mathcal{U}(\mathbf{M}_k(A))$ so that:

$$\rho_0(u) = u \otimes \mathbf{1}_k, \quad \text{and} \quad \rho_1(u) = u^k \oplus \mathbf{1}_{k-1}, \tag{1.6}$$

for any $u \in \mathcal{U}(A)$. This will be our starting point in the next section.

The following analogue of Theorem 3 should be a known result (conceivably, one might adapt the proof of Theorem 4.3 of [15] to this situation). In any case, we claim no originality for the result, and sketch a proof only to provide a guide to the next section.

**Theorem 1.7.** *Let $A$ be an approximately divisible $C^*$-algebra. Then $\pi_i(\mathcal{U}(A)) = K_{i-1}(A)$ for all integer $i \geq 0$.*

*Proof.* For any unital $C^*$-algebra $B$ and any integer $k \geq 1$, we define a map $\mu_k : \mathcal{U}(B) \to \mathcal{U}(\mathbf{M}_k(B))$ by letting $\mu_k(u) = u \oplus \mathbf{1}_{k-1}$ for any $u \in \mathcal{U}(B)$.

Let $A_n$, $F_n$ and $\theta_n$ be chosen as in Remark 1.4 to satisfy (1.2). For each $n$, we define a map $\eta_n : \mathcal{U}(\mathbf{M}_2(A_n)) \to \mathcal{U}(A_n \otimes F_n)$ as follows:

$$\eta_n(v) = \oplus_i [v^{k_i} \oplus \mathbf{1}_{k_i - 2}], \quad \forall v \in \mathcal{U}(\mathbf{M}_2(A_n)). \tag{1.7}$$

(This map was also used, for example, in [15].)

It is easy to check that $\eta_n$ is a well-defined map (cf. Remark 1.6). Moreover, note that the following diagram commutes up to homotopy:

$$\begin{array}{ccc} \mathcal{U}(A_n) & \xrightarrow{\cdot \otimes 1_{F_n}} & \mathcal{U}(A_n \otimes F_n) \\ \downarrow \mu_2 & \eta_n \nearrow & \downarrow \mu_2 \\ \mathcal{U}(\mathbf{M}_2(A_n)) & \xrightarrow{\cdot \otimes 1_{F_n}} & \mathcal{U}(\mathbf{M}_2(A_n \otimes F_n)). \end{array} \tag{1.8}$$

This follows easily from (1.6) (applied to $A$ and $\mathbf{M}_2(A)$).

Compare (1.8) with (1.2). Composing $\eta_n$ with the morphism $\theta_n$, and abusing the notation, we get a map $\eta_n : \mathcal{U}(\mathbf{M}_2(A_n)) \to \mathcal{U}(A_{n+1})$, and the diagram (1.8) still commutes up to homotopy when the algebra $A_n \otimes F_n$ (resp. the map $\cdot \otimes F_n$) there is replaced by $A_{n+1}$ (resp. the inclusion map).

Taking the limits (see also Lemma 2.7 below), we conclude that $\mu : \mathcal{U}(A) \to \mathcal{U}(\mathbf{M}_2(A))$ is a weak homotopy equivalence. Theorem 1.7 then follows easily. $\square$

2. Homotopy type of the unitary group

This section is devoted to the proof of Theorems 3. We follow closely the outline given in the proof of Theorem 1.7.

As we indicated there, the key to our proof is the construction of a map

$$\eta : \mathcal{U}(\mathbf{M}_k(A)) \to \mathcal{U}(A \otimes \mathcal{Z}), \tag{2.1}$$

with nice properties (cf. (1.8) in §1 and Corollary 2.6 below). Roughly speaking, $\eta$ is a continuous analogue of (1.7), and will be constructed in a way parallel to the constructions presented in (1.3) and (1.4). Note that by (1.5), the map $u \mapsto \rho(u)$ is in fact a continuous map from $\mathcal{U}(A)$ into $\mathcal{U}(A \otimes Z_{1,2})$. Clearly, this map is homotopic to the canonical map that sends $u$ to $u \otimes 1_{Z_{1,2}}$: A necessary homotopy can be given as follows:

$$h_{(s,t)}(u) = \rho_{st}(u); \quad (s,t) \in [0,1]^2. \tag{2.2}$$

We now generalize and juxtapose these constructions. Inspired by (1.3) and (1.4), we introduce the following

**Definition 2.1.** Let $n \geq 1$ be an integer. Suppose that $(W_1, \cdots, W_n)$ is a sequence of $n$ pathes in $\mathbf{SU}_n$, and for each $1 \leq j \leq n$,

$$W_j(0) = \mathbf{1}_n \quad \text{and} \quad W_j^*(1)\mathbf{e}_{j,j}W_j(1) = \mathbf{e}_{1,1}. \tag{2.3}$$

Then for any unital $C^*$-algebra, we define a map $\mathcal{W} : \mathcal{U}(A) \to \mathcal{U}(A \otimes Z_{1,n})$ by letting

$$\mathcal{W}(u;t) = \prod_{j=1}^{n} \left( W_j^*(t)[\mathbf{1}_n + (u-1)\mathbf{e}_{j,j}]W_j(t) \right), \tag{2.4}$$

for any $u \in \mathcal{U}(A)$ and $t \in [0,1]$. Such a map will be called *elementary*.

A few remarks are in order. As usual, $\mathbf{SU}_n$ in Definition 2.1 denotes the special unitary group of $\mathbf{M}_n$, and a path $W$ in $\mathbf{SU}_n$ means a *continuous* map $W : [0,1] \to \mathbf{SU}_n$. Paths that satisfy (2.3) are direct generalizations of the one defined in (1.3), and clearly, they always exist.

Also, recall that $\mathbf{M}_n$ can be identified canonically with the subalgebra $1_A \otimes \mathbf{M}_n$ of $\mathbf{M}_n(A)$ (cf. §1.1). This identification is used in (2.4) (and in (1.4)). It follows from (2.3) that, for any $u \in \mathcal{U}(A)$,

$$\mathcal{W}(u;0) = u \otimes \mathbf{1}_n \quad \text{and} \quad \mathcal{W}(u;1) = u^n \oplus \mathbf{1}_{n-1}, \tag{2.5}$$

and hence, $\mathcal{W}(u) \in \mathcal{U}(A \otimes Z_{1,n})$. It is easy to check that $\mathcal{W}$ given by (2.4) is a well-defined map (that is, it is continuous and basepoint-preserving). Incidentally, note that

$$\begin{aligned} W_j^*(t)[\mathbf{1}_n + (u-1)\mathbf{e}_{j,j}]W_j(t) \\ = 1_A \otimes \mathbf{1}_n + (u - 1_A) \otimes [W_j^*(t)e_{j,j}W_j(t)], \end{aligned} \tag{2.6}$$

for all $u \in \mathcal{U}(A)$, $t \in [0,1]$ and $1 \leq j \leq n$. This observation will be useful later (cf. the proof of Corollary 2.6 and (2.11) below).

Given $A$ and $n$, it might be tempting to work with one particular elementary map (as in (1.4)). But we shall need the flexibility offered by our slightly general approach (see the proof of Proposition 2.5). On the other hand, as we shall see now, there is essentially only one elementary map for the given $A$ and $n$.

**Lemma 2.2.** *Let $\mathcal{W}$, $\mathcal{V} : \mathcal{U}(A) \to \mathcal{U}(A \otimes Z_{1,n})$ be two elementary maps. Then there exists a continuous path of elementary maps $\mathcal{H}_t : \mathcal{U}(A) \to \mathcal{U}(A \otimes Z_{1,n})$ with $\mathcal{H}_0 = \mathcal{W}$ and $\mathcal{H}_1 = \mathcal{V}$.*

*Proof.* To fix notation, let $\mathcal{W}$ and $\mathcal{V}$ be given by two sequences $(W_1, \cdots, W_n)$ and $(V_1, \cdots, V_n)$, respectively. We shall show that for each $j$, $W_j$ is homotopic to $V_j$ in a suitable way.

For each integer $1 \leq j \leq n$, define a subset of $\mathbf{SU}_n$ as follows:
$$\mathbf{H}_{n,j} = \{w \in \mathbf{SU}_n : w^* \mathbf{e}_{j,j} w = \mathbf{e}_{1,1}\}.$$
Clearly, each $\mathbf{H}_{n,j}$ is a non-empty: $\mathbf{1}_n \in \mathbf{H}_{n,1}$, and $\mathbf{1}_n - (\mathbf{e}_{1,1} + \mathbf{e}_{j,j}) + (\mathbf{e}_{1,j} - \mathbf{e}_{j,1}) \in \mathbf{H}_{n,j}$ when $j > 1$ (cf. $W_1$ in (1.3)).

We claim that each $\mathbf{H}_{n,j}$ is connected. Indeed, one checks easily that
$$\mathbf{H}_{n,1} = \{ \begin{bmatrix} \det(w^*) & 0 \\ 0 & w \end{bmatrix} : w \in \mathbf{U}_{n-1}\},$$
where $\mathbf{U}_{n-1}$ is the unitary group. Since $\mathbf{U}_{n-1}$ is connected, so is $\mathbf{H}_{n,1}$. On the other hand, $\mathbf{H}_{n,j} = v \cdot \mathbf{H}_{n,1}$ for any $v \in \mathbf{H}_{n,j}$. This verifies the claim.

For each $j$, we now construct a homotopy between $W_j$ and $V_j$. Note that, by (2.3), both $W_j(1)$ and $V_j(1)$ are in $\mathbf{H}_{n,j}$. By the connectedness of $\mathbf{H}_{n,j}$, there is a path $h_j$ in $\mathbf{H}_{n,j}$ that joins $V_j(1)$ to $W_j(1)$ (that is, $h_j(0) = W_j(1)$, and $h_j(1) = V_j(1)$). Then, since $\mathbf{SU}_n$ is simply connected, there exists a continuous map $H_j : [0,1]^2 \to \mathbf{SU}_n$ such that
$$H_j(0, t) = W_j(t), \quad H_j(1, t) = V_j(t),$$
and
$$H_j(s, 0) = \mathbf{1}_n, \quad H_j(s, 1) = h_j(s) \in \mathbf{H}_{n,j},$$
for all $s$, $t \in [0, 1]$. In other words, $H_j$ is a homotopy between $W_j$ and $V_j$ (with paths that satisfy (2.3)).

Our conclusion follows immediately. □

We are now ready to construct (a local version of) the map $\eta$ promised in (2.1). In the construction, we use the following convention: If $f \in \mathcal{U}(\mathbf{M}_m(A) \otimes Z_{1,n})$, $g \in \mathcal{U}(\mathbf{M}_n(A) \otimes Z_{1,m})$, and $f(1) = g(1)$, then $f * g$ denotes the element in $\mathcal{U}(A \otimes Z_{m,n})$ given by:
$$(f * g)(t) = \begin{cases} f(2t), & \text{if } t \in [0, 1/2]; \\ g(2 - 2t), & \text{if } t \in [1/2, 1]. \end{cases} \tag{2.7}$$

More generally, let $\beta^{(0)} : \mathcal{U}(A) \to \mathcal{U}(\mathbf{M}_m(A) \otimes Z_{1,n})$ and $\beta^{(1)} : \mathcal{U}(A) \to \mathcal{U}(\mathbf{M}_n(A) \otimes Z_{1,m})$ be two maps with $\beta^{(0)}(u; 1) = \beta^{(1)}(u; 1)$ for all $u \in \mathcal{U}(A)$, then we define a map $\beta^{(0)} * \beta^{(1)} : \mathcal{U}(A) \to \mathcal{U}(A \otimes Z_{m,n})$ by letting:
$$(\beta^{(0)} * \beta^{(1)})(u) = \beta^{(0)}(u) * \beta^{(1)}(u). \tag{2.8}$$

**Definition 2.3.** Let $A$ be a unital $C^*$-algebra, $k \geq 1$ an integer, and $Z_{m,n}$ a prime dimension drop algebra with $m \geq k$ and $n \geq k$. Choose two elementary maps $\mathcal{W}_0 : \mathcal{U}(\mathbf{M}_m(A)) \to \mathcal{U}(\mathbf{M}_m(A) \otimes Z_{1,n})$ and $\mathcal{W}_1 : \mathcal{U}(\mathbf{M}_n(A)) \to \mathcal{U}(\mathbf{M}_n(A) \otimes Z_{1,m})$, and let
$$\eta(u) = \mathcal{W}_0(u^m \oplus \mathbf{1}_{(m-k)}) * \mathcal{W}_1(u^n \oplus \mathbf{1}_{(n-k)}),$$

for any $u \in \mathcal{U}(\mathbf{M}_k(A))$. It is easy to check that $\eta : \mathcal{U}(\mathbf{M}_k(A)) \to \mathcal{U}(A \otimes Z_{m,n})$ is a well-defined map (that is, $\eta$ is continuous and basepoint-preserving). Such a map will be called a *basic* map.

If $m = k = 1$, then any basic map $\eta : \mathcal{U}(A) \to \mathcal{U}(A \otimes Z_{1,n})$ is an elementary map. This follows from the fact that, for any unital $C^*$-algebra $B$, the only elementary map $\mathcal{W} : \mathcal{U}(B) \to \mathcal{U}(B \otimes Z_{1,1})$ is the identity map.

The auxiliary elementary maps $\mathcal{W}_j$ used in Definition 2.3 certainly exist and, by Lemma 2.2, are unique up to homotopy. Therefore, up to homotopy, the basic map $\eta$ in Definition 2.3 depends only on $A$, $k$, and $m$ and $n$.

Furthermore, as one might expect (cf. (2.2)), when $k = 1$, the map $\eta$ constructed in Definition 2.3 is homotopic to the natural embedding map $\cdot \otimes 1_{Z_{m,n}} : \mathcal{U}(A) \to \mathcal{U}(A \otimes Z_{m,n})$. This is the content of Proposition 2.5 below. To streamline the proof, we establish the following lemma.

**Lemma 2.4.** *Let $\mathcal{V} : A \to A \otimes Z_{1,m}$ and $\mathcal{W} : \mathbf{M}_m(A) \to \mathbf{M}_m(A) \otimes Z_{1,n}$ be two elementary maps. Define a map $\gamma : \mathcal{U}(A) \to \mathcal{U}(\mathbf{M}_m(A) \otimes Z_{1,n})$ by letting*

$$\gamma(u;t) = \mathcal{W}\big(\mathcal{V}(u;t);t\big), \quad \forall u \in \mathcal{U}(A).$$

*(1) There is a continuous family of maps $\Gamma_s : \mathcal{U}(A) \to \mathcal{U}(\mathbf{M}_m(A) \otimes Z_{1,n})$ such that $\Gamma_0 = \gamma$, $\Gamma_1(u) = \mathcal{W}(u^m \oplus \mathbf{1}_{m-1})$, and $\Gamma_s(u;1) = u^{mn} \oplus \mathbf{1}_{mn-1}$ for any $u \in \mathcal{U}(A)$ and $s \in [0,1]$.*

*(2) For any $u \in \mathcal{U}(A)$, $\gamma(u) \in \mathcal{U}(A \otimes Z_{1,mn})$. In fact, $\gamma : \mathcal{U}(A) \to \mathcal{U}(A \otimes Z_{1,mn})$ is an elementary map.*

*Proof.* Clearly, $\gamma : \mathcal{U}(A) \to \mathcal{U}(\mathbf{M}_m(A) \otimes Z_{1,n})$ is continuous, and $\gamma(1) = 1$.

The homotopy in part (1) can be constructed as follows:

$$\Gamma_s(u;t) = \mathcal{W}\big(\mathcal{V}(u; s + t - st); t\big),$$

for all $u \in \mathcal{U}(A)$ and $(s,t) \in [0,1]^2$. Direct computations show that the family $\Gamma_s$ satisfies all requirements in part (1) (note that $\mathcal{V}(u;1) = u^m \oplus \mathbf{1}_{m-1}$ and $\mathcal{W}(\mathcal{V}(u;1);1) = u^{mn} \oplus \mathbf{1}_{mn-1}$ for all $u \in \mathcal{U}(A)$). This proves part (1).

To prove part (2), note first that:

$$\gamma(u;0) = u \otimes \mathbf{1}_{mn}, \quad \forall u \in \mathcal{U}(A).$$

Therefore, $\gamma(u) \in \mathcal{U}(A \otimes Z_{1,mn})$. Furthermore, suppose that $\mathcal{V}$ and $\mathcal{W}$ are induced by the sequences $(V_1, \cdots, V_m)$ (of paths in $\mathbf{SU}_m$) and $(W_1, \cdots, W_n)$ (of paths in $\mathbf{SU}_n$), respectively (cf. Definition 2.1). Then a direct computation shows that $\gamma$ is induced by the sequence $V_i \otimes W_j$ (in the lexicographical order of $(i,j)$), where $V_i \otimes W_j$ denotes the path in $\mathbf{SU}_{mn}$ given by:

$$(V_i \otimes W_j)(t) = V_i(t) \otimes W_j(t), \quad t \in [0,1].$$

We leave the details to the reader. $\square$

We are now ready to prove the main technical result of this section:

**Proposition 2.5.** *Let $A$ be a unital $C^*$-algebra, and $Z_{m,n}$ a prime dimension drop algebra. Let $\iota : \mathcal{U}(A) \to \mathcal{U}(A \otimes Z_{m,n})$ be the natural embedding map given by $\iota(u) = u \otimes 1_{Z_{m,n}}$. Suppose that $\eta : \mathcal{U}(A) \to \mathcal{U}(A \otimes Z_{m,n})$ is a basic map (cf. Definition 2.3, where $k = 1$). Then $\eta$ and $\iota$ are homotopic (as maps from $\mathcal{U}(A)$ to $\mathcal{U}(A \otimes Z_{m,n})$).*

*Proof.* Let $\mathcal{W}_j$ be the auxiliary elementary maps used to construct $\eta$ (cf. Definition 2.3). Let $\mathcal{V}_0 : \mathcal{U}(A) \to \mathcal{U}(\mathbf{M}_m(A))$ and $\mathcal{V}_1 : \mathcal{U}(A) \to \mathcal{U}(\mathbf{M}_n(A))$ be two elementary maps. By Definition 2.3 (with $k = 1$) and Lemma 2.2, we have

$$\eta(u) = \mathcal{W}_0\big(\mathcal{V}_0(u; 1)\big) * \mathcal{W}_1\big(\mathcal{V}_1(u; 1)\big), \quad u \in \mathcal{U}(A).$$

Hence, by Lemma 2.6, there are two elementary maps $\gamma_0$ and $\gamma_1 : \mathcal{U}(A) \to \mathcal{U}(A \otimes Z_{1,mn})$ such that $\eta$ is homotopic to $\gamma_0 * \gamma_1$ (cf. (2.8)). By Lemma 2.2, we might as well assume that $\gamma_0 = \gamma_1$, which we denote simply by $\gamma$ from now on.

On the other hand, $\gamma * \gamma$ and $\iota$ are homotopic. To see this, define a continuous family of maps $\gamma_s : \mathcal{U}(A) \to \mathcal{U}(A \otimes Z_{1,mn})$ by $\gamma_s(u; t) = \gamma(u; st)$, for any $u \in \mathcal{U}(A)$ and $s, t \in [0, 1]$. Let $h_s = \gamma_s * \gamma_s$ (cf. (2.8)). One checks easily that $h$ is a homotopy between $\iota$ and $\gamma * \gamma$.

Therefore, $\eta$ and $\iota$ are homotopic. □

To state the next result, we make a few comments on the notation. For any unital $C^*$-algebra $A$, we denote by $\iota$ the natural unital embedding of $A$ into $A \otimes B$ (that is, $\iota(a) = a \otimes 1_B$), where $B$ is any unital nuclear $C^*$-algebra. Secondly, we use the same notation for a unital $*$-homomorphism and its induced map between the unitary groups. Finally, recall that for any unital $C^*$-algebra $A$, the map $\mu_k : \mathcal{U}(A) \to \mathcal{U}(\mathbf{M}_k(A))$ is given by

$$\mu_k(u) = u \oplus \mathbf{1}_{k-1}, \quad \forall u \in \mathcal{U}(A). \tag{2.9}$$

**Corollary 2.6.** *Let $A$ be a unital $C^*$-algebra and $k \geq 1$ an integer. Then there exists a map $\eta : \mathcal{U}(\mathbf{M}_k(A)) \to \mathcal{U}(A \otimes \mathcal{Z})$ such that the following diagram commutes up to homotopy:*

$$\begin{array}{ccc} \mathcal{U}(A) & \xrightarrow{\iota} & \mathcal{U}(A \otimes \mathcal{Z}) \\ \downarrow \mu_k & \eta \nearrow & \downarrow \mu_k \\ \mathcal{U}(\mathbf{M}_k(A)) & \xrightarrow{\iota \otimes \mathrm{id}_k} & \mathcal{U}(\mathbf{M}_k(A \otimes \mathcal{Z})), \end{array} \tag{2.10}$$

*where $\mathrm{id}_k$ is the identity map on $\mathbf{M}_k$.*

*Proof.* Let $Z_{m,n}$ be a prime dimension drop algebra with $m \geq k$ and $n \geq k$. Construct a basic map $\eta : \mathcal{U}(\mathbf{M}_k(A)) \to \mathcal{U}(A \otimes Z_{m,n})$ as in Definition 2.3. Clearly, by Proposition 1.2, it suffices to verify that the diagram (2.10) commutes up to homotopy when the algebra $\mathcal{Z}$ there is replaced by $Z_{m,n}$.

**The upper-left triangle.** Consider the two maps $\iota$ and $\eta \circ \mu_k : \mathcal{U}(A) \to \mathcal{U}(A \otimes Z_{m,n})$. A moment's reflection shows that $\eta \circ \mu_k$ is a basic map for which Proposition 2.5 applies. Therefore, $\eta \circ \mu_k$ is homotopic to $\iota$.

**The lower-right triangle.** The homotopy between $\mu_k \circ \eta$ and $\iota \otimes \mathrm{id}_k$ is also a consequence of Proposition 2.5 (applied to $B = \mathbf{M}_k(A)$). To see this, we first identify $\mathbf{M}_k(A \otimes Z_{m,n})$ with $B \otimes Z_{m,n}$ in a natural way: Note that (cf. §1.1)

$$\mathbf{M}_k(A \otimes Z_{m,n}) = (A \otimes Z_{m,n}) \otimes \mathbf{M}_k;$$

and

$$B \otimes Z_{m,n} = (A \otimes \mathbf{M}_k) \otimes Z_{m,n}.$$

We identify these two algebras by extending the following map:
$$(a \otimes f) \otimes c \mapsto (a \otimes c) \otimes f,$$
for $a \in A$, $f \in Z_{m,n}$, and $c \in \mathbf{M}_k$.

Under this identification, the map $\iota \otimes \mathrm{id}_k : B \to B \otimes Z_{m,n}$ is now given by $\iota$ (which maps $b$ to $b \otimes 1_{Z_{m,n}}$). We now calculate the composite map $\mu_k \circ \eta : \mathcal{U}(B) \to \mathcal{U}(B \otimes Z_{m,n})$, and show that it is homotopic to a basic map.

First we fix notation. Let $\mathcal{W}_0$ and $\mathcal{W}_1$ be the two auxiliary elementary maps used to construct $\eta$ (cf. Definition 2.3). Moreover, we assume that the map $\mathcal{W}_0$ is induced by the sequence $(W_1, \cdots, W_n)$ of $n$ paths in $\mathbf{SU}_n$ (cf. Definition 2.1). Clearly, in exactly the same way, this sequence also defines an elementary map $\mathcal{W}_0^B : \mathcal{U}(\mathbf{M}_m(B)) \to \mathcal{U}(\mathbf{M}_m(B) \otimes Z_{1,n})$. Similarly, we have also an elementary map $\mathcal{W}_1^B : \mathcal{U}(\mathbf{M}_n(B)) \to \mathcal{U}(\mathbf{M}_n(B) \otimes Z_{1,m})$.

Recall that $B = \mathbf{M}_k(A)$, hence $\mathbf{M}_k(B) = A \otimes \mathbf{M}_k \otimes \mathbf{M}_k$. Let $\psi$ be the inner automorphism on $\mathbf{M}_k(B)$ induced by the flip automorphism on $\mathbf{M}_k$, that is,
$$\psi(a \otimes e \otimes f) = a \otimes f \otimes e,$$
for all $a \in A$ and $e, f \in \mathbf{M}_k$. Finally, let $\Psi : \mathcal{U}(B) \to \mathcal{U}(\mathbf{M}_k(B))$ be the map given by $\Psi(v) = \psi(v \oplus \mathbf{1}_{k-1})$ for any $v \in \mathcal{U}(B)$.

We are now ready to have a formula for $\mu_k \circ \eta$. Note that $\mu_k$ is a point-wise operation that is multiplicative and $*$-preserving. A careful calculation, based on this fact and on (2.6), reveals that:
$$(\mu_k \circ \eta)(v) = \mathcal{W}_0^B\big(\Psi(v)^m \oplus \mathbf{1}_{(m-k)}\big) * \mathcal{W}_1^B\big(\Psi(v)^n \oplus \mathbf{1}_{(n-k)}\big),$$
for any $v \in \mathcal{U}(B)$.

On the other hand, since $\psi$ is induced by a flip automorphism on $\mathbf{M}_k$, it can be implemented by a unitary in $\mathbf{M}_k \otimes \mathbf{M}_k \subseteq \mathbf{M}_k(B)$. It follows that $\psi$ is homotopic (as an automorphism) to the identity map on $\mathbf{M}_k(B)$. Therefore, the two maps $\Psi$ and $\mu_k : \mathcal{U}(B) \to \mathcal{U}(\mathbf{M}_k(B))$ are homotopic. Clearly, this homotopy induces a homotopy between $\mu_k \circ \eta$ and the following map:
$$\mathcal{U}(B) \ni v \mapsto \mathcal{W}_0^B(v^m \oplus \mathbf{1}_{(m-1)}) * \mathcal{W}_1^B(v^n \oplus \mathbf{1}_{(n-1)}) \in \mathcal{U}(B \otimes Z_{m,n}).$$

Clearly, the latter is a basic map and, by Proposition 2.5, is homotopic to the natural map $\iota : \mathcal{U}(B) \to \mathcal{U}(B \otimes Z_{m,n})$. Therefore, $\mu_k \circ \eta$ and $\iota$ are homotopic. The proof is completed. □

The following result should be well-known (cf. the unitary analogues of Proposition 3.3.3, Proposition 3.4.3, and §9.2.3 of [1]).

**Lemma 2.7.** *Let $A$ be a unital $C^*$-algebra and $\{A_n\}$ an increasing sequence of unital $C^*$-subalgebras with $\cup_n A_n$ dense in $A$. Then the canonical map*
$$\mathrm{incl.}_* : \lim_{n \to +\infty} \pi_i(\mathcal{U}(A_n)) \to \pi_i(\mathcal{U}(A))$$
*is an isomorphism for each integer $i \geq 0$.*

*Proof.* The case $i = 0$ is standard. For $i > 0$, let $\infty \in S^k$ be the basepoint of the sphere $S^k$. It is easy to show that given any unitary $u \in C(S^k; A)$ with $u(\infty) = 1$, and any constant $\epsilon > 0$, there is a unitary $u_n \in C(S^k; A_n)$ with $u_n(\infty) = 1$ and $||u - u_n|| < 1$. The rest of the proof is again standard. □

**Theorem 2.8.** *Let $A$ be a unital $\mathcal{Z}$-stable $C^*$-algebra, $k > 0$ an integer. Then the canonical embedding map $\mu_k : \mathcal{U}(A) \to \mathcal{U}(\mathbf{M}_k(A))$ (cf. (2.9)) is a weak homotopy equivalence, that is, $\mu_k$ induces an isomorphism*

$$\mu_* : \pi_i(\mathcal{U}(A)) \cong \pi_i(\mathcal{U}(\mathbf{M}_k(A)))$$

*of homotopy groups for any integer $i \geq 0$.*

*Proof.* Since $A$ is $\mathcal{Z}$-stable, by Remark 1.4, there is an increasing sequence $\{A_n\}$ of unital $C^*$-subalgebras of $A$, with $\cup_n A_n$ dense in $A$, such that for each $n$, there is an isomorphism $\theta_n : A_n \otimes \mathcal{Z} \cong A_{n+1}$ that makes diagram (1.1) commute (that is, $a = \theta_n(a \otimes 1_\mathcal{Z})$ for all $a \in A_n$).

By Corollary 2.6, there are maps $\eta_n : \mathcal{U}(\mathbf{M}_k(A_n)) \to \mathcal{U}(A_{n+1})$ such that the following diagram commutes up to homotopy:

$$\begin{array}{ccccccc}
\mathcal{U}(A_1) & \xrightarrow{\text{incl.}} & \mathcal{U}(A_2) & \xrightarrow{\text{incl.}} & \mathcal{U}(A_3) & \xrightarrow{\text{incl.}} & \cdots \\
\downarrow \mu_k & \eta_1 \nearrow & \downarrow \mu_k & \eta_2 \nearrow & \downarrow \mu_k & \nearrow & \cdots \\
\mathcal{U}(\mathbf{M}_k(A_1)) & \xrightarrow{\text{incl.}} & \mathcal{U}(\mathbf{M}_k(A_2)) & \xrightarrow{\text{incl.}} & \mathcal{U}(\mathbf{M}_k(A_2)) & \xrightarrow{\text{incl.}} & \cdots
\end{array}$$

It follows from Lemma 2.7 that, for each integer $i \geq 0$, the sequence $(\mu_k)_*$ has a limit $\mu_* : \pi_i(\mathcal{U}(\mathbf{M}_k(A))) \to \pi_i(\mathcal{U}(A))$, and $\mu_*$ and $(\mu_k)_*$ are inverse to each other. This completes the proof. $\square$

Theorems 2 and 3 follow immediately from Theorem 2.8 and Bott periodicity.

*Remark 2.9.* Based on Theorem 2.8, one can show that $\mu_k$ is actually a homotopy equivalence. See the proof of Lemma 4.1 of [17] for details.

It is known that, for any $C^*$-algebra $A$, $K_1(A) \cong K(A \otimes \mathcal{Z})$ (cf. Lemma 2.11 of [11]). Therefore, from Theorem 2, we have:

**Corollary 2.10.** *For any unital $C^*$-algebra $A$, $K_1(A) \cong \mathcal{U}(A \otimes \mathcal{Z})/\mathcal{U}_0(A \otimes \mathcal{Z})$.*

*Remark 2.11.* For a sample list of earlier works similar to Theorem 3, we mention Cuntz [7], Rieffel[12], Thomsen[15], and Zhang [17], [18]. See [15] and [17] for nice surveys of this topic.

It is known ([11]) that $\mathcal{Z}$-stable $C^*$-algebras include all unital simple separable infinite dimensional AF-algebras and all unital separable simple nuclear purely infinite $C^*$-algebras. It is also clear that if $A$ is $\mathcal{Z}$-stable, then $A \otimes B$ is also $\mathcal{Z}$-stable for any $C^*$-algebra $B$. Therefore, the main result of this section, Theorem 3, covers many cases treated earlier.

In connection with the results in [17] and [18], we mention the following open question:

**Question 2.12.** *Let $A$ be a unital separable simple nuclear infinite dimensional $C^*$-algebra of real rank zero.*
 1. *Is $A$ $\mathcal{Z}$-stable?*
 2. *Is $A$ approximately divisible?*

The second part of the question was first raised in [4].

In the rest of this section, we consider an analogue of Theorem 3 for non-unital $C^*$-algebras. For this purpose, we first extend the definition of $\mathcal{U}(A)$.

Let $A$ be a non-unital $C^*$-algebra. We denote by $A^+$ the unitization of $A$ (cf. §3.2 of [1]), and define
$$\mathcal{U}(A) = \mathcal{U}(A^+) \cap (1_{A^+} + A).$$

Clearly, $\mathcal{U}(A)$ is homeomorphic to the group $U(A)$ of quasi-unitaries in $A$ (cf. [15]). Again, we consider $\mathcal{U}(A)$ as a topological space with basepoint $1_{A^+}$. This definition is in fact functorial. In particular, if $\phi : A \to B$ is a $C^*$-homomorphism between two non-unital $C^*$-algebras, then there is a canonical (continuous basepoint-preserving) map $\mathcal{U}(\phi) : \mathcal{U}(A) \to \mathcal{U}(B)$, given by
$$\mathcal{U}(\phi)(1_{A^+} + a) = 1_{B^+} + \phi(a) \quad \text{if } a \in A \text{ and } 1_{A^+} + a \in \mathcal{U}(A).$$

In other words, $\mathcal{U}(\phi)$ is the restriction of the canonical map $\phi^+ : A^+ \to B^+$. For brevity, we denote $\mathcal{U}(\phi)$ simply by $\phi$.

With this definition, we claim that Theorem 3 extends to the non-unital case, with essentially the same proof. We discuss briefly the necessary changes.

Our proof of Theorem 2.8 (hence of Theorem 3) has three components: Corollary 2.6, Lemma 2.7 and a structure result for $\mathcal{Z}$-stable $C^*$-algebras (Remark 1.4). Note that the structure result is a consequence of a theorem on $\mathcal{Z}$ (Theorem 1.3), and is true for all $\mathcal{Z}$-stable $C^*$-algebras.

It is also easy to check Lemma 2.7 for a non-unital $C^*$-algebra $A$. In fact, the proof becomes easier. The case $i = 0$ again follows from functional calculus (note that two close unitaries in $\mathcal{U}(A)$ are homotopic within $\mathcal{U}(A)$). The case $i > 0$ can be reduced to the case $i = 0$, with a standard trick:
$$\pi_i(\mathcal{U}(A)) = \pi_0(\mathcal{U}(C_0(\mathbb{R}^i \otimes A)).$$

It is slightly more involved to verify Corollary 2.6 for a non-unital $C^*$-algebra $A$. But the idea is simple and natural, namely, unitization and restriction. We elaborate.

First, we check that maps in the following diagram (cf. (2.10)) are well-defined:
$$\begin{array}{ccc} \mathcal{U}(A) & \xrightarrow{\iota} & \mathcal{U}(A \otimes \mathcal{Z}) \\ \downarrow \mu_k & & \downarrow \mu_k \\ \mathcal{U}(\mathbf{M}_k(A)) & \xrightarrow{\iota \otimes \mathrm{id}_k} & \mathcal{U}(\mathbf{M}_k(A \otimes \mathcal{Z})). \end{array}$$

The horizontal maps are induced by $*$-homomorphisms and, as we discussed earlier, they are well-defined. For the vertical maps, we start with the map $\mu_k : \mathcal{U}(A^+) \to \mathcal{U}(\mathbf{M}_k(A^+))$ (cf. (2.9)). It is easy to see that $\mu_k(\mathcal{U}(A)) \subseteq \mathcal{U}(\mathbf{M}_k(A))$. Hence, by restriction, we have a map $\mu_k : \mathcal{U}(A) \to \mathcal{U}(\mathbf{M}_k(A))$. To get the other vertical map, apply the same procedure to the non-unital algebra $A \otimes \mathcal{Z}$. Incidentally, this diagram commutes.

Corollary 2.6 (the non-unital case) says that there exists a map $\eta : \mathcal{U}(\mathbf{M}_k(A)) \to \mathcal{U}(A \otimes \mathcal{Z})$ such that the diagram (2.10) commutes up to homotopy. To construct this map, we follow exactly the same route as in the unital case. Here it is important to note that, if $\mathcal{W} : \mathcal{U}(A^+) \to \mathcal{U}(A^+ \otimes Z_{1,n})$ is any elementary map, then
$$\mathcal{W}(\mathcal{U}(A)) \subseteq \mathcal{U}(A \otimes Z_{1,n}). \tag{2.11}$$

This follows immediately from Definition 2.1 (see also (2.6)). Therefore, again by restriction, we have a map $\mathcal{W} : \mathcal{U}(A) \to \mathcal{U}(A \otimes Z_{1,n})$ (recall that a map here means a continuous and basepoint-preserving map). We then construct a basic map $\eta$ as in the unital case.

The rest of the proof follows the same pattern. Namely, we repeat the proof for the unital case, and check that every map there (including those in a homotopy) has, through unitization and restriction, a counterpart in the non-unital case. This process is straightforward, and not particularly illuminating. We leave the details to the interested reader.

In conclusion, Theorem 3 also holds for non-unital $\mathcal{Z}$-stable $C^*$-algebras.

## 3. Comparison of full projections

Let $A$ be a $C^*$-algebra. Recall that a projection $p \in A$ is called full, if $p$ is not contained in any proper closed two-sided ideal of $A$. In this section, we apply Theorem 2 to study comparison theory for full projections in $\mathcal{Z}$-stable $C^*$-algebras, and, in particular, we prove Theorem 1.

**Notation 3.1.** First we fix the notation. Let $A$ be a unital $\mathcal{Z}$-stable $C^*$-algebra. Let $A_n$ be the increasing sequence of unital $*$-subalgebras of $A$ constructed in Remark 1.4. Recall that
$$A = \overline{\cup_n A_n}. \tag{3.1}$$
See condition (1) in Remark 1.4. Also note that, for each $n$, there is an isomorphism $\theta_n : A_n \otimes \mathcal{Z} \cong A$ such that the following diagram commutes:

$$\begin{array}{ccc} A_n & \xrightarrow{\cdot \otimes 1_{\mathcal{Z}}} & A_n \otimes \mathcal{Z} \\ \downarrow \text{incl.} & & \downarrow \theta_n \\ A & = & A. \end{array} \tag{3.2}$$

(Compare with condition (2) in Remark 1.4.) It follows that the inclusion $A_n \hookrightarrow A$ induces an isomorphism on K-theory (cf. Lemma 2.11 of [11]). In particular, if $p$ and $q$ are two projections in $A_n$, then $[p] = [q]$ in $K_0(A)$ if and only if $[p] = [q]$ in $K_0(A_n)$.

We now prove two basic lemmas.

**Lemma 3.2.** Let $A$ be any $\mathcal{Z}$-stable $C^*$-algebra. If $p \in A$ is a projection, then $pAp$ is also $\mathcal{Z}$-stable.

*Proof.* We use the notation in §3.1. Suppose that $p \in A$ is a projection. By (3.1), and by functional calculus (cf. Proposition 4.5.2 of [1]), we can find a projection $q \in A_n$ for some $n$ large enough, such that $||p - q|| < 1$. Hence $p$ and $q$ are unitarily equivalent, and $pAp$ and $qAq$ are isomorphic. Therefore, without loss of generality, we may assume that $p \in A_n$. Then, by (3.2),
$$pAp \cong (pA_n p) \otimes \mathcal{Z}.$$
This completes the proof. $\square$

*Remark 3.3.* Let $A$ and $B$ be two *unital* $C^*$-algebras that are strong Morita equivalent to each other. Then $A$ is $\mathcal{Z}$-stable if and only if $B$ is $\mathcal{Z}$-stable. This follows immediately from Lemma 3.2, and the obvious fact that $\mathbf{M}_n(A) = A \otimes \mathbf{M}_n$ is $\mathcal{Z}$-stable for any $n \geq 1$ if $A$ is $\mathcal{Z}$-stable.

The situation is less clear for general (non-unital) $C^*$-algebras.

The next lemma is a local version of Theorem 1(a). As in [1], here we use $\sim$ (resp. $\sim_u$, and $\sim_h$) to denote Murray-von Neumann equivalence (resp. unitary equivalence, and homotopy) of projections. We shall also use freely the following elementary fact: Let $p$ and $q$ be two projections in $A$. If $(p \oplus \mathbf{0}_k) \sim (q \oplus \mathbf{0}_k)$ in $\mathbf{M}_{k+1}(A)$ for some integer $k \geq 0$, then $p \sim q$ in $A$.

**Lemma 3.4.** *Let $A$ be a unital $\mathcal{Z}$-stable $C^*$-algebra, and $p$, $q$ two full projections in $A$. If $[p] = [q] \in K_0(A)$, then $p \otimes 1_{\mathcal{Z}} \sim q \otimes 1_{\mathcal{Z}}$ in $A \otimes \mathcal{Z}$.*

*Proof.* By a Theorem of Blackadar (cf. Theorem 3.1.4 of [2]), there exist an integer $k_0$ such that $p \otimes 1_k \sim q \otimes 1_k$ (in $\mathbf{M}_k(A)$) for all $k > k_0$. Replacing $A$ by $\mathbf{M}_2(A)$ and $p$ (resp. $q$) by $p \oplus 0$ (resp. $q \oplus 0$) if necessary, we assume that $p \otimes 1_k \sim_h q \otimes 1_k$ (in $\mathbf{M}_k(A)$) for all $k > k_0$.

Let $m$ and $n$ be a pair of relatively prime integers with $m > k_0$ and $n > k_0$. It follows from the preceding paragraph that there exist unitaries $u_0 \in \mathbf{M}_m(A)$ and $u_1 \in \mathbf{M}_n(A)$ such that
$$p \otimes \mathbf{1}_m = u_0^*(q \otimes \mathbf{1}_m)u_0, \quad \text{and} \quad p \otimes \mathbf{1}_n = u_1^*(q \otimes \mathbf{1}_n)u_1. \tag{3.3}$$
Our next objective is to find a continuous path $U_t$ in $\mathcal{U}(\mathbf{M}_{mn}(A))$ connecting $u_1 \otimes \mathbf{1}_m$ to $u_0 \otimes \mathbf{1}_n$ such that
$$p \otimes \mathbf{1}_{mn} = U_t^*(q \otimes \mathbf{1}_{mn})U_t, \quad \forall t \in [0,1]. \tag{3.4}$$

Unfortunately, this is not always possible. For example, a necessary (but not sufficient) condition is that $n[u_0] = m[u_1]$ in $K_1(A)$. We first modify $u_0$ and $u_1$, as follows.

Let $V = (u_1 \otimes \mathbf{1}_m)^* \cdot (u_0 \otimes \mathbf{1}_n)$. Then, by (3.3), we have
$$V^*(p \otimes \mathbf{1}_{mn})V = p \otimes \mathbf{1}_{mn}. \tag{3.5}$$
Let $A_p = pAp$, and $V_p = (p \otimes \mathbf{1}_{mn})V(p \otimes \mathbf{1}_{mn})$. By (3.5), $V_p$ is a unitary in $\mathbf{M}_{mn}(A_p)$. It follows from Lemma 3.2 and Theorem 2 that there is a unitary $w_p$ in $A_p$ such that $[w_p] = -[V_p] \in K_1(A_p)$. Similarly, there is a unitary $w_{1-p}$ in $A_{1-p} = (1-p)A(1-p)$ such that $[w_{1-p}] = -[(1-p) \otimes \mathbf{1}_{mn} \cdot V \cdot (1-p) \otimes \mathbf{1}_{mn}] \in K_1(A_{1-p})$. On the other hand, since $m$ and $n$ are relatively prime, we can choose (and fix) a pair $(j, k)$ of integers such that
$$j \cdot m + k \cdot n = 1.$$
Finally, set:
$$\tilde{u}_0 = u_0 \cdot \big((w_p + w_{1-p})^k \oplus \mathbf{1}_{m-1}\big] \in \mathcal{U}(\mathbf{M}_m(A)),$$
and
$$\tilde{u}_1 = u_1 \cdot [(w_e + w_{1-e})^{-j} \oplus \mathbf{1}n - 1] \in \mathcal{U}(\mathbf{M}_n(A)).$$
Note that (3.3) still holds with $u_0$ (resp. $u_1$) there replaced by $\tilde{u}_0$ (resp. $\tilde{u}_1$).

We claim that there is a continuous path $U_t$ in $\mathcal{U}(\mathbf{M}_{mn}(A))$ such that:
$$U_0 = \tilde{u}_0 \otimes \mathbf{1}_n, \quad U_1 = \tilde{u}_1 \otimes \mathbf{1}_m, \tag{3.6}$$
and $U_t$ satisfies (3.4) for each $t$.

To verify the claim, let $\widetilde{V} = (\tilde{u}_1 \otimes \mathbf{1}_m)^* \cdot (\tilde{u}_0 \otimes \mathbf{1}_n)$. Then a careful bookkeeping shows that:
$$[(p \otimes \mathbf{1}_{mn}) \cdot \widetilde{V} \cdot (p \otimes \mathbf{1}_{mn})] = m \cdot j \cdot [w_p] + [V_p] + n \cdot l \cdot [w_p] = 0 \in K_1(A_p).$$

Again, by Lemma 3.2 and Theorem 2, there is a continuous path $V_p(t)$ (of unitaries in $\mathbf{M}_{mn}(A_p)$) such that

$$V_p(0) = (p \otimes \mathbf{1}_{mn}) \cdot \widetilde{V} \cdot (p \otimes \mathbf{1}_{mn}), \quad \text{and} \quad V_p(1) = 1 \in \mathbf{M}_{mn}(A_p).$$

Similarly, there is a continuous path $V_{1-p}(t)$ (of unitaries in $\mathbf{M}_{mn}(A_{1-p})$) that connects the identity $1 \in \mathbf{M}_{mn}(A_{1-p})$ to $(1 - p \otimes \mathbf{1}_{mn}) \cdot \widetilde{V} \cdot (1 - p \otimes \mathbf{1}_{mn})$. Now let

$$U_t = (\tilde{u}_1 \otimes \mathbf{1}_n) \cdot (V_p(t) + V_{1-p}(t)), \quad \forall t \in [0,1].$$

It is easy to verify that $U_t$ satisfies (3.6) and (3.4). In other words, $U$ defines a unitary in $A \otimes Z_{m,n}$, and

$$p \otimes 1_{Z_{m,n}} = U^*(q \otimes 1_{Z_{m,n}})U.$$

By Proposition 1.2, we have $p \otimes 1_{\mathcal{Z}} \sim_u q \otimes 1_{\mathcal{Z}}$. In particular, $p \otimes 1_{\mathcal{Z}} \sim q \otimes 1_{\mathcal{Z}}$. $\square$

Note that in the first step of the proof, we might have to pass to $\mathbf{M}_2(A)$. Therefore, we can not conclude that $p \otimes 1_{\mathcal{Z}} \sim_u q \otimes 1_{\mathcal{Z}}$, which, in fact, is not true in general (for example, it fails in infinite $C^*$-algebras).

Let $A$ be as in Lemma 3.4 and $k \geq 1$ is any integer. Then $\mathbf{M}_k(A)$ is unital and $\mathcal{Z}$-stable. Therefore, the conclusion of Lemma 3.4 also holds for full projections in a matrix algebra over $A$. Incidentally, in Lemma 3.4, the $\mathcal{Z}$-stability condition on $A$ can be dropped. Note that, by Theorem 1.3, $A \otimes \mathcal{Z}$ is always $\mathcal{Z}$-stable, and $A \otimes \mathcal{Z} \otimes \mathcal{Z} \cong A \otimes \mathcal{Z}$. Hence, if $A$ is not $\mathcal{Z}$-stable, we pass to $A \otimes \mathcal{Z}$ and apply Lemma 3.4 to get the same conclusion.

The following result was stated as Theorem 1(a) in the introduction.

**Theorem 3.5.** *Let $A$ be a unital $\mathcal{Z}$-stable $C^*$-algebra. If $p$ and $q$ are two full projections in $A$, and $[p] = [q]$ in $K_0(A)$, then $p \sim q$.*

*Proof.* We use the notation in §3.1. As in the proof of Lemma 3.2, we assume, without loss of generality, that $p$ and $q \in A_n$, for some $n$ large enough. From the discussion in §3.1, it follows that $[p] = [q] \in K_0(A_n)$. Hence, by (3.2) and Lemma 3.4, $p \sim q$ in $A$. This completes the proof. $\square$

Before turning to Theorem 1(b), we discuss some ramifications of Theorem 3.5.

**Corollary 3.6.** *Let $A$ be a unital $\mathcal{Z}$-stable $C^*$-algebra and $v \in A$ a partial isometry. If both $1 - v^*v$ and $1 - vv^*$ are full projections, then there is a partial isometry $v^\perp \in A$ such that $v + v^\perp \in \mathcal{U}_0(A)$.*

*Proof.* Let $p = v^*v$, $q = vv^*$. Clearly, $[p] = [q]$ and hence $[1-p] = [1-q]$ in $K_0(A)$. By assumption, $1 - p$ and $1 - q$ are full, hence Theorem 3.5 applies to ensure that $1 - p \sim 1 - q$. Let $w \in A$ be a partial isometry with

$$w^*w = 1 - p \quad \text{and} \quad ww^* = 1 - q.$$

Clearly, $v + w$ is an unitary in $A$.

Denote $1 - q$ by $q^\perp$. Since $q^\perp$ is a full projection, the inclusion $q^\perp A q^\perp \hookrightarrow A$ induces an isomorphism of K-theory (cf. proof of Corollary 2.7 in [6]). In particular, by Lemma 3.2 and Theorem 2, there is an unitary $u \in \mathcal{U}(q^\perp A q^\perp)$ such that $[\tilde{u}] = -[v+w] \in K_1(A)$, where $\tilde{u} = q + u \in A$. Let $v^\perp = \tilde{u}w$. Note that

$$[v + v^\perp] = [\tilde{u}(v+w)] = 0 \in K_1(A).$$

Therefore, by Theorem 2, $v + v^\perp \in \mathcal{U}_0(A)$. This completes the proof. $\square$

In the following result, $K_0(A)$ is considered as a (scaled) preordered group in the usual way (cf. §6.1 of [1]). Also, as usual, for two projections $p$ and $q$ in the same algebra, we write $p \leq q$ if $p$ is a subprojection of $q$, and write $p \preceq q$, if there is a projection $q_0$ such that $p \sim q_0 \leq q$.

**Corollary 3.7.** *Let $A$ be a unital $\mathcal{Z}$-stable $C^*$-algebra, and $p$, $q$ two full projections in $A$.*

(1) *If $1 - p$ and $1 - q$ are also full, and if $[p] = [q] \in K_0(A)$, then $p \sim_h q$.*
(2) *If $[p] \leq [q]$ (in $K_0(A)$), then $p \preceq q$.*

*Proof.* (1) This follows immediately from Theorem 3.5 and Corollary 3.6.

(2) Since $[p] \leq [q]$, there exists a projections $r \in \mathbf{M}_k(A)$ (for some $k \geq 1$) such that $[p] + [r] = [q]$. By Theorem 3.5, $p \oplus r \sim q \oplus \mathbf{0}_k$ in $\mathbf{M}_{k+1}(A)$. Hence, $p \preceq q$. This completes the proof. □

In particular, if $A$ is unital $\mathcal{Z}$-stable and *simple*, and $p$ and $q$ are two non-trivial projections (that is, $p, q \notin \{0, 1\}$) in $A$, then $[p] = [q] \in K_0(A)$ if and only if $p \sim_h q$. For similar earlier results, see Proposition 3.4 of [16].

This also brings us to the Grassmann space $\mathcal{P}(A)$ of $A$. Recall that $\mathcal{P}(A)$ is the space of nontrivial projections in $A$, endowed with the relative norm topology from $A$. Through the map $p \mapsto 2p - 1$, it is homeomorphic to the space of non-trivial symmetries in $A$. Note that this space could be empty. For example, $\mathcal{P}(\mathcal{Z}) = \emptyset$. However, when not empty, it carries interesting information. The preceding paragraph is a statement on the connected components of this space. We now consider higher homotopy groups. We take an arbitrary point in $\mathcal{P}(A)$ as its basepoint.

**Proposition 3.8.** *Let $A$ be a unital simple $\mathcal{Z}$-stable $C^*$-algebra. If $\mathcal{P}(A) \neq \emptyset$, then*

$$\pi_i(\mathcal{P}(A)) = \begin{cases} K_0(A) & \text{if } i \text{ is even,} \\ K_1(A) & \text{if } i \text{ is odd,} \end{cases} \quad \forall i \geq 1.$$

*Proof.* The proof uses Theorem 3, and follows the same line given in Sections 3.11, 3.12, and 3.13 of [17]. □

In connection with Proposition 3.8, it might be of interest to note that $\mathbf{M}_2(A)$ obviously contains non-trivial projection if $A$ is unital (and, as usual, non-trivial).

To prove Theorem 1.(b), We recall some standard notation (cf. [5]). Let $A$ be a unital $C^*$-algebra. We denote by $\mathbf{M}_\infty(A)$ the union $\cup_n \mathbf{M}_n(A)$, where the inclusion $\mathbf{M}_n(A) \hookrightarrow \mathbf{M}_{n+1}(A)$ is given by $a \mapsto a \oplus 0$. Let $V(A)$ denote the scaled preordered semigroup of all equivalence class of projections in $\mathbf{M}_\infty(A)$. For a projection $p$, we denote by $[p]$ its class in $V(A)$.

For discussions on quasi-traces and related matters, see [3], [14], [4] and [5]. Note also that if $A$ is exact, any quasi-trace on $A$ is actually a trace (cf. [10]).

The following is a more precise statement of Theorem 1(b) in the introduction.

**Theorem 3.9.** *Let $A$ be a unital $\mathcal{Z}$-stable $C^*$-algebra, and $p$ and $q$ two projections in $A$ with $q$ full. If $\tau(p) < \tau(q)$ for all quasi-traces $\tau$ on $A$, then $p \preceq q$.*

*Proof.* As in the proofs of Lemma 3.2 and Theorem 3.5(1), by invoking (3.1) and (3.2), we only need to show that $p \otimes 1_\mathcal{Z} \preceq q \otimes 1_\mathcal{Z}$ in $A \otimes \mathcal{Z}$. This will be similar to the proof of Lemma 3.4.

Note that $[q]$ is an order unit in $V(A)$, since $q$ is full. Invoking successively three results (Theorem 3.3, Lemma 2.8, and Lemma 2.3) in [5], we conclude that $k \cdot [p] + [1] \leq k \cdot [q]$ in $V(A)$ for some integer $k > 0$. Clearly, $[p] \leq [1]$, and therefore, $(k+1) \cdot [p] \leq k \cdot [q]$. Using again the fact that $[q]$ is an order unit, we can easily find an integer $k_0 \geq 1$ such that $(k+1) \cdot [p] \leq k \cdot [q]$ for all integer $k \geq k_0$ (see, for example, [2]).

Let $m$ and $n$ be a pair of relatively prime integers strictly larger than $k_0$. Then, in particular, $m \cdot [p] \leq (m-1) \cdot [q]$. That is, there is a partial isometry $v_0 \in \mathbf{M}_m(A)$ such that:

$$v_0^* v_0 = p \otimes \mathbf{1}_m, \quad \text{and} \quad v_0 v_0^* \leq (q \otimes \mathbf{1}_{m-1}) \oplus 0. \tag{3.7}$$

Similarly, there is a partial isometry $v_1 \in \mathbf{M}_n(A)$ such that:

$$v_1^* v_1 = p \otimes \mathbf{1}_n, \quad \text{and} \quad v_1 v_1^* \leq (q \otimes \mathbf{1}_{n-1}) \oplus 0. \tag{3.8}$$

Let $W = (v_1 \otimes \mathbf{1}_m) \cdot (v_0 \otimes \mathbf{1}_n)^* \in \mathbf{M}_{mn}(A)$.

By (3.7) and (3.8), $W$ is actually a partial isometry in $B$, where $B = \mathbf{M}_{mn}(qAq)$. Clearly, $B$ is unital and, by Lemma 3.2, $\mathcal{Z}$-stable. Using (3.7) and (3.8) again, one can check easily that $W$ satisfies the conditions in Corollary 3.6. Therefore, there is a unitary $U \in \mathcal{U}_0(B)$ such that $W = U \cdot (W^*W)$. Let $U_t$ be a continuous path of unitaries in $B$ that connects $U$ to $1_B$.

Since $B$ is a corner of $\mathbf{M}_{mn}(A)$, we now regard $U_t$ as a continuous path of partial isometries in $\mathbf{M}_{mn}(A)$. Let $V_t = U_t \cdot (v_0 \otimes \mathbf{1}_n)$ for $t \in [0, 1]$. By (3.7), $V_t$ is a continuous path of partial isometries. By the choice of $U_t$, we have:

$$V_0 = v_0 \otimes \mathbf{1}_n, \quad \text{and} \quad V_1 = v_1 \otimes \mathbf{1}_m.$$

In other words, $V$ defines a partial isometry in $A \otimes Z_{m,n}$. A direct calculation shows that:

$$V^* V = p \otimes 1_{Z_{m,n}}, \quad \text{and} \quad V V^* \leq q \otimes 1_{Z_{m,n}}.$$

(The latter follows since $U_t \in B$ for each $t$.) It follows that $p \otimes 1_{Z_{m,n}} \preceq q \otimes 1_{Z_{m,n}}$. By Proposition 1.2, we have $p \otimes 1_{\mathcal{Z}} \preceq q \otimes 1_{\mathcal{Z}}$. This completes the proof. $\square$

In particular, if $A$ is unital $\mathcal{Z}$-stable and admits no quasi-traces, then $1_A \preceq q$ for any full projection $q$. This also follows from Corollary 3.7.

*Remark 3.10.* Our proof of Theorem 3.9 refines the proof of Theorem 1 of [9]. In fact, when $A$ is *simple*, then Theorem 3.9 follows from Theorem 1 of [9], Corollary 3.7 above, and some results in [14]. We sketch the proof. First we note that if $A$ is not stably finite, then by Theorem 3 of [9], $A$ is purely infinite. In this case, the conclusion of Theorem 3.9 is well-known (cf. [7]). So we assume that $A$ is stably finite. Then by Theorem 1 of [9] and Corollary 3.7(2) above, $V(A)$ is almost unperforated in the sense of Rordam (cf. [14]). The conclusion then follows from Theorem 6.1 and Proposition 3.2 of [14].

*Remark 3.11.* As we mentioned in Remark 3.4, Theorem 3.5 and Theorem 3.9 also hold for full projections in a matrix algebra over a unital $\mathcal{Z}$-stable $C^*$-algebra.

We conclude this note with a brief discussion on analogous results for non-unital $C^*$-algebras.

In Theorem 3.5, the hypothesis that $A$ be unital is redundant. If $A$ is $\mathcal{Z}$-stable, and if $p$ and $q$ are two full projections in $A$ with $[p] = [q] \in K_0(A)$, then it is easy to check that the two projections $(p \oplus 0)$ and $(0 \oplus q)$ in $(p \oplus q) \mathbf{M}_2(A) (p \oplus q)$ satisfy all conditions in Theorem 3.5, hence they are equivalent. It follows easily that $p \sim q$.

Theorem 3.9 also has an analogue for general $\mathcal{Z}$-stable algebras. Namely, if $A$ is a $\mathcal{Z}$-stable $C^*$-algebra, if $p$ and $q$ are two projections in $A$ with $q$ full, and if $\tau(p) < \tau(q)$ for every (non-zero) quasitrace $\tau$ on the Pedersen ideal of $A$, then $p \preceq q$. This follows from the trick used in the preceding paragraph, and an extension result due to Blackadar and Handelman (Proposition II.4.2 of [3]).

# References

[1] B. Blackadar, *K-Theory for Operator Algebras*, MSRI Publ. Vol. 5, Springer-Verlag, New York, 1985.
[2] B. Blackadar, *Rational $C^*$-algebras and nonstable K-theory*, Rocky Mountain Math. J. 20(1990), 285–316.
[3] B. Blackadar and D. Handelman, *Dimension functions and traces on $C^*$-algebras*, J. Funct. Anal., 45 (1982), 297–340.
[4] B. Blackadar, A. Kumjian and M. Rordam, *Approximately central matrix units and the structure of noncommutative tori*, K-theory, 6(1992), 267–284.
[5] B. Blackadar and M. Rordam, *Extending states on preordered semigroups and the existence of quasitraces on $C^*$-algebras*, J. Algebra, 152 (1992), 240–247.
[6] L. Brown, *Stable isomorphism of hereditary subalgebras of $C^*$-algebras*, Pacific J. Math., 71(1977), 335–348.
[7] J. Cuntz, *K-theory for certain $C^*$-algebras*, Ann. of Math. (2) 113(1981), 181–197.
[8] G. A. Elliott, *The classification problem for amenable $C^*$-algebras*, Proceedings of the International Congress of Mathematicians, Vol. 2 (Zurich, 1994), 922–932, Birkhauser, Basel, 1995.
[9] G. Gong, X. Jiang and H. Su, *Obstructions to $\mathcal{Z}$-stability for unital simple $C^*$-algebras*, preprint.
[10] U. Haagerup, *Quasitraces on exact $C^*$-algebras are traces*, to appear.
[11] X. Jiang and H. Su, *On a simple unital projectionless $C^*$-algebra*, preprint.
[12] M. A. Rieffel, *The homotopy groups of the unitary groups of non-commutative tori*, J. Operator Theory, 17(1987), 237–254.
[13] M. Rordam, *On the structure of simple $C^*$-algebras tensored with a $UHF$-algebra I*, J. Funct. Anal. 100(1991), 1–17.
[14] M. Rordam, *On the structure of simple $C^*$-algebras tensored with a $UHF$-algebra II*, J. Funct. Anal. 107(1992), 255–269.
[15] K. Thomsen, *Non-stable K-theory for operator algebras*, K-theory, 4(1991), 245–267.
[16] S. Zhang, *Diagonalizing projections in the multiplier algebras and matrices over a $C^*$-algebra*, Pacific J. Math., 145(1990), 181–200.
[17] S. Zhang, *Matricial structure and homotopy type of simple $C^*$-algebras with real rank zero*, J. Operator Theory 26 (1991), 283–312.
[18] S. Zhang, *On the homotopy type of the unitary group anf the Grassmann space of purely infinite simple $C^*$-algebras*, K-theory, to appear (accepted in 1992).

The Fields Institute, 222 College Street, Toronto, Ontario, Canada M5T 3J1
*E-mail address*: `jiang@fields.utoronto.ca`